\documentclass{amsart}

\usepackage{amsmath,amssymb,amsthm}

\newtheorem{prop}{Proposition}
\newtheorem{thm}[prop]{Theorem}

\theoremstyle{definition}

\newtheorem*{rem}{Remark}

\def\co{\colon\thinspace}

\newcommand{\N}{\mathbb N}

\newcommand{\R}{\mathbb R}

\newcommand{\Z}{\mathbb Z}

\DeclareMathOperator{\Cont}{Cont}

\DeclareMathOperator{\Diff}{Diff}

\DeclareMathOperator{\st}{st}

\newcommand{\xist}{\xi_{\st}}

\begin{document}

\author{Hansj\"org Geiges}
\author{Mirko Klukas}
\address{Mathematisches Institut, Universit\"at zu K\"oln,
Weyertal 86--90, 50931 K\"oln, Germany}
\email{geiges@math.uni-koeln.de, mklukas@math.uni-koeln.de}

\title[On the space of contact structures]{The fundamental group of
the space of contact structures on the $3$-torus}

\date{}

\begin{abstract}
We show that the fundamental group of the space of contact structures
on the $3$-torus (based at the standard contact structure)
is isomorphic to the integers.
\end{abstract}


\maketitle


\section{Introduction}
On the $3$-torus $T^3$ with coordinates $x,y,\theta\in\R/2\pi\Z$
the standard tight contact structure is given by
\[ \xist=\ker (\cos\theta \, dx-\sin\theta\, dy).\]
In~\cite{gego04} it was shown that the space $\Xi (T^3)$
of contact structures on $T^3$ has nontrivial topology:
its fundamental group based at $\xist$ contains a subgroup
isomorphic to the integers~$\Z$.
An alternative proof of that result was given by
Bourgeois~\cite{bour06}. He also found nontrivial
subgroups in higher homotopy groups of certain higher-dimensional
contact manifolds. These subgroups are detected by their
nontrivial action on the contact homology of the manifold
in question.

Here we give a quick proof that $\Z$ constitutes the full fundamental group.

\begin{thm}
\label{thm:main}
The fundamental group of $\Xi (T^3)$, based at~$\xist$,
is isomorphic to~$\Z$.
\end{thm}

The same statement holds with $\xist$ replaced by
\[ \xi_n=\ker(\cos(n\theta)\, dx-\sin (n\theta)\, dy)\]
for any $n\in\N$.

Our proof of this theorem is based on a homotopy exact sequence
that will be described in Section~\ref{section:exact}. A similar
approach was used in~\cite{dige10} to show that the fundamental group
of the space of contact structures on $S^1\times S^2$, based
at the standard contact structure $\ker (z\, d\theta + x\, dy-y\, dx)$,
is isomorphic to~$\Z$.

Eliashberg~\cite{elia91} has shown
that the space of tight contact structures on the $3$-sphere
fixed at a point is contractible; this will in fact be a key
ingredient in the proof of Theorem~\ref{thm:main}.
Beyond these results, very little is known about comparable
homotopical questions.

In 2005, Bourgeois announced a proof of Theorem~\ref{thm:main}
that employs quite different methods: instead of relating the homotopy
type of $\Xi (T^3)$ to that of the diffeomorphism resp.\ contactomorphism
group of $(T^3,\xist)$ via a homotopy exact sequence, a direct
approach based on a subtle cut-and-paste procedure along
surfaces is used to show that $\pi_1 (\Xi (T^3),\xist)$ is
contained in~$\Z$. We learned only recently that details of this proof
have been worked out by Fabien Ng\^o in his MSc thesis~\cite{ngo05}.
See also the comment following Proposition~\ref{prop:Cont}.

\section{A homotopy exact sequence}
\label{section:exact}
Let $(M,\xi_0)$ be a closed contact manifold. Write $\Xi_0$
for the component of the space of contact structures
on~$M$ containing~$\xi_0$. By $\Diff_0$ we denote the identity
component of the diffeomorphism group of $M$, and by
$\Cont_0$ its subgroup of contactomorphisms of~$\xi_0$. (Beware that
$\Cont_0$ is not, in general, connected.)
As shown in~\cite{gego04}, the inclusion map
$i\co\Cont_0\rightarrow \Diff_0$
and the map
\[ \begin{array}{rcl}
\sigma\co \Diff_0 & \longrightarrow & \Xi_0\\
\phi & \longmapsto & \phi_*\xi_0,
\end{array} \]
induce a homotopy exact sequence
\[... \stackrel{\Delta}{\longrightarrow} \pi_k(\Cont_0)
\stackrel{i_{\#}}{\longrightarrow} \pi_k(\Diff_0)
\stackrel{\sigma_{\#}}{\longrightarrow} \pi_k(\Xi_0)
\stackrel{\Delta}{\longrightarrow} \pi_{k-1}(\Cont_0)
\stackrel{i_{\#}}{\longrightarrow} ...\]

For $M=T^3$ we have $\pi_1(\Diff_0)\cong\Z^3$, generated by the
full turns around the three $S^1$-factors, cf.~\cite{hatc80}.
The shifts in the $x$- and $y$-direction are contactomorphisms of~$\xist$,
so the subgroup $\Z^2+\{0\}\subset\Z^3$ lies in the image of
\[ i_{\#}\co \pi_1(\Cont_0)
\longrightarrow \pi_1(\Diff_0),\]
and hence in the kernel of
\[ \sigma_{\#}\co \pi_1(\Diff_0)
\longrightarrow \pi_1(\Xi_0).\]
So the relevant part of the homotopy exact sequence for
$(T^3,\xi_0)$ reduces to
\[ \Z\stackrel{\sigma_{\#}}{\longrightarrow}
\pi_1(\Xi_0)\stackrel{\Delta}{\longrightarrow}\pi_0(\Cont_0),\]
where $\Z=\{(0,0)\}+\Z\subset\Z^3$ is generated by the loop
\[ (x,y,\theta)\mapsto (x,y,\theta+2\pi s),\;\; s\in [0,1].\]
As shown in~\cite{gego04}, $\sigma_{\#}(\Z)$ is an infinite cyclic subgroup
of $\pi_1(\Xi_0)$.
\section{Proof of Theorem~\ref{thm:main}}
Theorem~\ref{thm:main} is now an immediate consequence
of the following result.

\begin{prop}
\label{prop:Cont}
For $(T^3,\xist)$, the space $\Cont_0$ is connected, i.e.\
every contactomorphism of $(T^3,\xist)$ that is topologically
isotopic to the identity is also isotopic to the
identity via contactomorphisms.
\end{prop}

This proposition, in turn, is a direct corollary of a result
of Giroux~\cite[Th\'eo\-r\`eme~4]{giro01}. A proof of
Theorem~\ref{thm:main} based on this line of reasoning was
given in~\cite{gego}. Unfortunately, the published proof
of Giroux's result is incomplete. (The error occurs
in \cite[Proposition~10]{giro01}. The proofs of the main
results, though, can be fixed using the methods of
Patrick Massot's thesis~\cite{mass08}). Here we present
a proof of Proposition~\ref{prop:Cont} based on a result of
Ghiggini~\cite{ghig06}. For the language of convex surfaces
in the sense of Giroux see~\cite{giro91} or~\cite{geig08}.

\begin{proof}[Proof of Proposition~\ref{prop:Cont}]
The $2$-torus $T:=\{y=0\}$ in~$T^3$ is a vertical torus in standard
form in the sense of~\cite{ghig06}: a convex torus (transverse
to the contact vector field~$\partial_y$) with two dividing curves given
by $\theta\in\{ 0,\pi\}$ and linear Legendrian ruling given
by~$\partial_{\theta}$ (the requirement `vertical' means that
the Legendrian ruling has a nontrivial $\partial_{\theta}$-component).
Likewise, the $2$-torus $S:=\{x=0\}$ is such a vertical torus.

Now let $\phi\in\Cont_0$ be given. Then $\phi(T)$ is smoothly isotopic
to~$T$, and by \cite[Lemma~6.5]{ghig06} it is contact isotopic to~$T$.
So $\phi$ is isotopic in $\Cont_0$ to a contactomorphism $\phi'$
fixing~$T$. By the isotopy extension theorem for surfaces
in contact $3$-manifolds, cf.~\cite[Theorem~2.6.13]{geig08},
we may in fact assume that $\phi'$ fixes $T$ pointwise.

The second vertical torus $S$ intersects $T$ in a Legendrian
ruling curve, and it is smoothly isotopic to the torus
$\phi'(S)$ sharing this property. By \cite[Lemma~6.3]{ghig06},
$\phi'(S)$ is contact isotopic to~$S$ (keeping
$T$ fixed). Thus, we can isotope $\phi'$ in $\Cont_0$ to
a contactomorphism $\phi''$ fixing both $T$ and~$S$ pointwise.

The solid torus $T^3\setminus (S\cup T)$ becomes,
after edge-rounding~\cite[Lemma~3.11]{hond00}, a standard tight solid
torus with two dividing curves of slope~$-1$.
Now consider a convex meridional disc $D$ in this solid torus and
its image under~$\phi''$. Note that the dividing set of $D$
consists of a single properly embedded arc. By isotopy
discretisation~\cite[Lemma~3.10]{hond02}, an idea going back to
Colin, one sees that $\phi''(D)$ is contact isotopic to~$D$,
so we can now isotope $\phi''$ in $\Cont_0$ to $\phi'''$
fixing $S,T$ and~$D$. (This isotopy discretisation is also
the basis for Ghiggini's results.)

The complement $T^3\setminus (S\cup T\cup D)$, after edge-rounding,
is a standard tight $3$-ball~$(B^3,\xist)$, and $\phi'''$
may be regarded as a contactomorphism of $(B^3,\xist)$
fixing the boundary. According to Eliashberg~\cite[Theorem~2.4.2]{elia91},
the space of tight contact structures on $B^3$ with fixed
boundary condition is contractible,
in particular, its fundamental group is trivial.
A homotopy exact sequence as in the previous section then 
shows that the space of contactomorphisms of $(B^3,\xist)$
fixed at the boundary (and topologically isotopic to the identity)
is connected. This concludes the proof of the proposition.
\end{proof}

\begin{rem}
In the last part of the foregoing proof we only had to deal
with contactomorphisms of $B^3$ (fixed at the boundary)
that are known by assumption to be topologically isotopic
to the identity. In fact, however, by Hatcher's proof of the
Smale conjecture (saying that the diffeomorphism group of $S^3$
retracts to the orthogonal group of isometries) this topological
condition is always satisfied: the full diffeomorphism group
of $B^3$ relative to the boundary is contractible. It follows that
the same is true for the contactomorphism group of $(B^3,\xist)$
relative to the boundary, cf.~\cite[Th\'eor\`eme~17]{giro01}.
\end{rem}

\end{document}